\documentclass[psamsfont,a4paper,12pt]{amsart}
\usepackage[english]{babel}
\numberwithin{equation}{section}

\usepackage{setspace}
\usepackage{enumerate}
\usepackage{afterpage}
\usepackage[usenames,dvipsnames]{color}
\usepackage{graphicx}

\usepackage{amssymb,latexsym,amsmath}
\usepackage{amssymb}
\usepackage{mathrsfs}
\usepackage{amsfonts}
\usepackage{latexsym}
\usepackage[latin1]{inputenc}
\usepackage{pstricks,multido,pst-plot}
\usepackage[none]{hyphenat}

\usepackage{amsthm}
\usepackage{stackrel}                         

\usepackage{tabularx}
\usepackage{longtable}

\usepackage{verbatim}
\usepackage{blindtext}
\usepackage{multicol,fancyheadings}
\usepackage{float}
\usepackage[vcentermath,enableskew]{youngtab}
\usepackage{subfig}

\allowdisplaybreaks

\newtheorem{theorem}{Theorem}[section]			

\newtheorem{Thm}{Theorem}[]					

\newtheorem*{thm}{Theorem}					

\theoremstyle{definition}

\theoremstyle{remark}

\newtheorem{Proof}[Thm]{Proof}

\usepackage{tikz}
\usetikzlibrary{arrows}

\usepackage{paralist}
\usepackage{cite}

\title{Non-standard binary representations and the Stern sequence}

\author[Anders]{Katie Anders}
\address{Department of Mathematics, University of Texas at Tyler, Tyler, TX 75799}
\email{kanders@uttyler.edu}

\author[Dawsey]{Madeline Locus Dawsey}
\address{Department of Mathematics, University of Texas at Tyler, Tyler, TX 75799}
\email{mdawsey@uttyler.edu}

\author[Gupta]{Rajat Gupta}
\address{Department of Mathematics, University of Texas at Tyler, Tyler, TX 75799}
\email{rgupta@uttyler.edu}

\author[Vandehey]{Joseph Vandehey}
\address{Department of Mathematics, University of Texas at Tyler, Tyler, TX 75799}
\email{jvandehey@uttyler.edu}

\begin{document}


\begin{abstract}
We show that the number of short binary signed-digit representations of an integer $n$ is equal to the $n$-th term in the Stern sequence. Various proofs are provided, including direct, bijective, and generating function proofs. We also show that this result can be derived from recent work of Monroe on binary signed-digit representations of a fixed length.
\end{abstract}

\maketitle

\section{Introduction}\label{introduction}

Balanced base-$b$ expansions, which allow both positive and negative digits, were introduced in 1726 by Colson \cite{colson1726short}, where they were known as ``negativo-affirmative arithmetick." In this paper, we will focus on binary signed-digit representations. A \textit{binary signed-digit (BSD) representation} of $n$ is of the form
\begin{equation}
n= \sum_{j=0}^\infty \epsilon_j 2^j, \qquad \epsilon_j\in\{1,0,-1\}. \label{eq:initial BSD definition}
\end{equation}
BSD representations are useful in computation, where they are implemented in general arithmetic \cite{avizienis1961signed} and in cryptography \cite{okeya,ruan}. Ergodic properties of BSD representations have also been studied in \cite{dajani2020invariant,dajani2006ergodic}.

Several authors have noted a connection between various binary representations and the Stern sequence, which is given by $s(0)=0$, $s(1)=1$, $s(2n)=s(n)$, and $s(2n+1)=s(n+1)+s(n)$, and so begins
\[
	0, 1, 1, 2, 1, 3, 2, 3, 1, 4, 3, 5, 2, 5, \ldots.
\] 
Considering Equation \eqref{eq:initial BSD definition} but allowing $\epsilon_j\in\{0,1,2\}$ instead of $\epsilon_j\in\{-1,0,1\}$ gives a \textit{hyperbinary representation} of $n$. Reznick \cite{BR} showed that if $f_{\mathrm{HB}}(n)$ is the number of hyperbinary representations of $n$, then $f_{\mathrm{HB}}(n) = s(n+1)$. More recently, Monroe \cite{ML} considered $f_\mathrm{BSD}(n,i)$, the number of BSD representations of $n$ using exactly $i$ digits:
\[
n= \sum_{j=0}^{i-1} \epsilon_j 2^j, \qquad \epsilon_j\in\{1,0,-1\}. 
\]
Monroe showed that $f_\mathrm{BSD}(n,i)=s(2^i-n)$. In this paper, we will combine the perspectives of Monroe and Reznick and show a direct connection between a restricted class of BSD representations and the Stern sequence.

To ease the readability of BSD representations, we will use $T$ to represent $-1$. (Some places prefer $\overline{1}$.) Further, we will often represent BSD representations using just the digits, as we do with standard binary representations, so that, for example,
\[
[10T1]_2 = 1\cdot 2^3+0\cdot 2^2+(-1)\cdot 2^1+1\cdot 2^0.
\]
The term $\epsilon_j 2^j$, where $j$ is as large as possible with $\epsilon_j \neq 0$, is said to be the \textit{leading term} of this expansion.

There exist infinitely many BSD representations for any non-zero integer. For example, $5$ can be written as
\[
[101]_2=[1T01]_2=[1TT01]_2=[1TTT01]_2=\cdots,
\]
because $2^i-2^{i-1}-2^{i-2}-\cdots-1 = 1$. This is why Monroe considered only BSD representations with a finite number of digits. We will consider a different restriction. We will say a BSD representation is \emph{short} if its leading two terms are neither $1T$ nor $T1$. Let $\overline{f}_\mathrm{BSD}(n)$ denote the number of short BSD representations of $n$. Then we see that $\overline{f}_\mathrm{BSD}(5)=3$ as $[101]_2$, $[11T]_2$, and $[10TT]_2$ are the only short BSD representations of $5$, while $\overline{f}_\mathrm{BSD}(8)=1$ since $[1000]_2$ is the only short BSD representation of $8$. Note that $f_\mathrm{HB}(n)=0$ for all $n<0$, as it is impossible to represent a negative number as a sum of non-negative values; however, $\overline{f}_\mathrm{BSD}(n)=\overline{f}_\mathrm{BSD}(-n)$ for all $n\in\mathbb{Z}$, since one can obtain a short BSD representation of $-n$ from a short BSD representation of $n$ simply by replacing all $1$'s with $T$'s and vice versa.

Our main results are the following:
\begin{theorem}\label{main1}
For all integers $n\ge 0$, we have $\overline{f}_\mathrm{BSD}(n)=s(n)$.
\end{theorem}
\begin{theorem}\label{main2} For all integers $n\ge 0$, we have $\overline{f}_\mathrm{BSD}(n)=f_\mathrm{HB}(n-1)$. \end{theorem}

Each of these theorems follows from the other combined with Reznick's result that $f_{\mathrm{HB}}(n)=s(n+1)$. We include both Theorem \ref{main1} and Theorem \ref{main2} because we will be giving a variety of proofs of these facts from different perspectives, and sometimes one theorem is more natural than the other. We will provide a direct proof, in the style of Reznick, of Theorem \ref{main1} in Section \ref{section2} and three proofs of Theorem \ref{main2}, including both a bijective proof and a generating function proof, in Section \ref{section3}.

We will also connect $\overline{f}_\mathrm{BSD}(n)$ with the function $f(n,i)$ studied by Monroe.
\begin{theorem}\label{main3}
    Let $n\ge 1$ be an integer and let $i=\lfloor \log_2 n\rfloor+1$ denote the number of digits in the standard binary representation of $n$. Then 
    \[
    f_\mathrm{BSD}(n,i+1)-f_\mathrm{BSD}(n,i)=\overline{f}_\mathrm{BSD}(n).
    \]
\end{theorem}
When combined with a recurrence relation of the Stern sequence, this will provide yet another proof of Theorem \ref{main2}.  We prove Theorem \ref{main3} in Section \ref{section4}.

\section{Proof of Theorem \ref{main1}}\label{section2}

In \cite{BR}, Reznick showed that $f_\mathrm{HB}(n-1)=s(n)$ for any positive integer $n$.  We employ the same style of proof to prove Theorem \ref{main1}, which we restate here.

\begin{thm}\label{BSD Stern correspondence} 
For all integers $n\geq 0$, we have $\overline{f}_\mathrm{BSD}(n)=s(n)$.
\end{thm}

We first consider a few base cases.  Note that $\overline{f}_\mathrm{BSD}(0)=0, \overline{f}_\mathrm{BSD}(1)=1, \overline{f}_\mathrm{BSD}(2)=1, \overline{f}_\mathrm{BSD}(3)=2, \overline{f}_\mathrm{BSD}(4)=1$, and $\overline{f}_\mathrm{BSD}(5)=3$. Suppose that the result holds for every non-negative integer less than $n$. By our calculations above, we may assume that $n\ge 6$.  Since $\{0,1\}\subset \{0,1,T\}$ and every non-negative integer has a standard binary representation, we know $\overline{f}_\mathrm{BSD}(n)\neq 0$.  Write 
\[
n=\sum_{j=0}^{\infty}\epsilon_j2^j=\epsilon_0+\epsilon_12^1+\epsilon_22^2+\epsilon_32^3+\cdots,
\]
where each $\epsilon_j\in\{0,1,T\}$.

Suppose $n$ is even.  Then there exists a non-negative integer $k<n$ such that $n=2k$.  We must have $\epsilon_0=0$, and
\[
n=\epsilon_12^1+\epsilon_22^2+\epsilon_32^3+\cdots=2\left(\epsilon_1+\epsilon_22^1+\epsilon_32^2+\cdots\right)=2k.
\]
Since appending $0$ to the end of a short BSD representation will not change the fact that it is short, we get that $\overline{f}_\mathrm{BSD}(n)=\overline{f}_\mathrm{BSD}(k)$, and by the strong inductive hypothesis, there are $\overline{f}_\mathrm{BSD}(k)=s(k)$ ways to write $k$. 

Now suppose $n$ is odd.  Then there exists a non-negative integer $k<n$ such that $n=2k+1$.  Since $n$ is odd, either $\epsilon_0=1$ or $\epsilon_0=T$.  Suppose $\epsilon_0=1$.  We have
\[
n=2k+1=1+\epsilon_12^1+\epsilon_22^2+\epsilon_32^3+\cdots,
\]
so
\[
2k=2\left(\epsilon_1+\epsilon_22^1+\epsilon_32^2+\cdots\right),
\]
and there are $\overline{f}_\mathrm{BSD}(k)$ ways to write $k$.
Now suppose $\epsilon_0=T$.  Then we have
\[
n=2k+1=T+\epsilon_12^1+\epsilon_22^2+\epsilon_32^3+\cdots,
\]
so
\[
2k+2=\epsilon_12^1+\epsilon_22^2+\epsilon_32^3+\cdots,
\]
and
\[
2(k+1)=2\left(\epsilon_1+\epsilon_22^2+\epsilon_32^2+\cdots\right).
\]
Of course there are $\overline{f}_\mathrm{BSD}(k+1)$ ways to write $k+1$. Appending $1$ or $T$ to a short BSD representation will not change the fact that it is short unless we are appending $T$ to $[1]_2$ or appending $1$ to $[T]_2$. Since we have assumed $n\ge 6$, we have $k\ge 3$, and thus these exceptional cases do not occur here. Thus $\overline{f}_\mathrm{BSD}(n)=\overline{f}_\mathrm{BSD}(2k+1)=\overline{f}_\mathrm{BSD}(k)+\overline{f}_\mathrm{BSD}(k+1)=s(k)+s(k+1)$, where we have used the strong inductive hypothesis.

Using the recurrence relations for the Stern sequence, we see that when $n=2k$, we have $\overline{f}_\mathrm{BSD}(n)=\overline{f}_\mathrm{BSD}(k)=s(k)=s(n)$, and when $n=2k+1$, we have $\overline{f}_\mathrm{BSD}(n)=s(k)+s(k+1)=s(2k+1)=s(n)$.  In all cases, $\overline{f}_\mathrm{BSD}(n)=s(n)$.\qed

\section{Proofs of Theorem \ref{main2}}\label{section3}

We recall Theorem \ref{main2}, and then we offer three different proofs.  Proof \ref{indirect} is an indirect proof combining Theorem \ref{main1} of this paper and Theorem 5.2 of \cite{BR}.  Proof \ref{bijection} explicitly constructs a bijection.  Proof \ref{generating_function_proof} proves a generating function identity from which the result follows.

\begin{thm}\label{BSD HB correspondence} 
For all integers $n\geq0$, we have $\overline{f}_\mathrm{BSD}(n)=f_\mathrm{HB}(n-1)$.
\end{thm}

\begin{Proof}\label{indirect}
As mentioned above, Reznick proved in \cite{BR} that $f_\mathrm{HB}(n-1)=s(n)$ for any positive integer $n$.  Since we know from Theorem  \ref{main1} that $\overline{f}_\mathrm{BSD}(n)=s(n)$ for any non-negative integer $n$, we can conclude that $\overline{f}_\mathrm{BSD}(n)=f_\mathrm{HB}(n-1)$.\qed
\end{Proof}

We also prove Theorem \ref{main2} by constructing a bijection between the set of hyperbinary representations of any non-negative integer $n$ and the set of short BSD representations of $n+1$.

\begin{Proof}\label{bijection}
Let $n$ be a non-negative integer, and suppose we have a hyperbinary representation of $n$ of the form $$n=\sum_{i=0}^\ell\epsilon_i2^i,$$ where each $\epsilon_i\in\{0,1,2\}$ and $\epsilon_\ell\neq0$.  In particular, any $n$ of this form is a positive integer.  Recall that we denote this representation by $[\epsilon_\ell\;\epsilon_{\ell-1}\;\cdots\;\epsilon_2\;\epsilon_1\;\epsilon_0]_2$. Let $+_2$ denote component-wise addition of binary representations, and define a function $f$ by
\begin{align*}
f([\epsilon_\ell\;\epsilon_{\ell-1}\;\cdots\;\epsilon_2\;\epsilon_1\;\epsilon_0]_2)&=[\epsilon_\ell\;\epsilon_{\ell-1}\;\cdots\;\epsilon_2\;\epsilon_1\;\epsilon_0]_2+_2[1\;\underbrace{T\cdots\;T}_{\ell+1\text{ digits}}]_2\\
&=[1\;(\epsilon_\ell+T)\;(\epsilon_{\ell-1}+T)\;\cdots\;(\epsilon_2+T)\;(\epsilon_1+T)\;(\epsilon_0+T)]_2.
\end{align*}
 Also define $f([0]_2)=[1]_2$.
Then $f$ maps a hyperbinary representation of $n$ with $\ell+1$ digits to a short BSD representation of $n+1$ with $\ell+2$ digits, since $$[1\;\underbrace{T\;\cdots\;T}_{\ell+1\text{ digits}}]_2=2^{\ell+1}-\sum_{i=0}^\ell2^i=2^{\ell+1}-(2^{\ell+1}-1)=1,$$ each $\epsilon_i+T\in\{0,1,T\}$, and $\epsilon_\ell+T\neq T$.   We will now show that $f$ is a bijection between the set of all hyperbinary representations of $n$, for $n\geq 0$, and the set of all short BSD representations of $n+1$, for $n\ge 0$.  To see that $f$ is onto, consider a short BSD representation of $n+1$.  Such a representation must be of the form $[1\; \delta_{\ell-1} \; \delta_{\ell-2} \; \cdots \; \delta_1 \; \delta_0]_2$ with $\ell\ge 0$, where the leading term is $1$ since $n+1$ is positive. In the case where $\ell=0$, the only short BSD representation of this form is $[1]_2$, and $f([0]_2)=[1]_2$ as noted above. In the case where $\ell\geq 1$, $[(\delta_{\ell-1}+1)\;(\delta_{\ell-2}+1)\;\cdots\;(\delta_1+1)\;(\delta_0+1)]_2$ is a hyperbinary representation of $n$ since $\delta_i+1\in \{0,1,2\}$ for each $i$, and we have
\[
f([(\delta_{\ell-1}+1)\;(\delta_{\ell-2}+1)\;\cdots\;(\delta_1+1)\;(\delta_0+1)]_2) = [1\; \delta_{\ell-1} \; \delta_{\ell-2} \; \cdots \; \delta_1 \; \delta_0]_2.
\]
Thus $f$ is onto.   It is also clear that $f$ is injective, since the only way that the image of two hyperbinary representations $n=[\epsilon_\ell\;\cdots\;\epsilon_0]$ and $n=[\epsilon_\ell'\;\cdots\;\epsilon_0']$ can be equal is if each $\epsilon_i=\epsilon_i'$.  Thus $f$ is a bijection, and this proves that $\overline{f}_{\mathrm{BSD}}(n+1)=f_{\mathrm{HB}}(n)$ for any non-negative integer $n$.\qed
\end{Proof}

Alternatively, we can construct a map from the set of short BSD representations of $n+1$ to the set of hyperbinary representations of $n$.  Such a map $g$ satisfying $g=f^{-1}$ is defined by
\begin{align*}
g([\delta_k\;\delta_{k-1}\;\cdots\;\delta_2\;\delta_1\;\delta_0]_2)&=[\delta_k\;\delta_{k-1}\;\cdots\;\delta_2\;\delta_1\;\delta_0]_2+_2[T\;\underbrace{1\cdots\;1}_{k\text{ digits}}]_2\\
&=[(\delta_{k-1}+1)\;\cdots\;(\delta_2+1)\;(\delta_1+1)\;(\delta_0+1)]_2
\end{align*}
and $g([1]_2):=[0]_2$.  The fact that $g$ is a bijection also proves Theorem \ref{main2}.

Our third proof of Theorem \ref{main2} employs a generating function argument.

\begin{Proof}\label{generating_function_proof}
Here we prove Theorem \ref{main2} by proving the following formal generating function identity:
\begin{equation}\label{generating_function_identity}
q\prod_{n=0}^\infty \left( 1+ q^{2^n}+q^{2\cdot 2^n}\right) = q+\sum_{N=1}^\infty \left[\prod_{i=0}^{N-2} \left(1+q^{2^i}+q^{-2^i}\right)\right] \left(1+q^{2^{N-1}}\right)q^{2^N}.
\end{equation}
Note that these are the desired generating functions for hyperbinary representations and short BSD representations. The series on the left side of \eqref{generating_function_identity} is the generating function of $f_\mathrm{HB}(n-1)$.  This is well known and can be found in the proof of Theorem 3.1 in \cite{Northshield}, among other places. Observe that the $n-1$ is reflected in the extra copy of $q$ out front.  

On the other hand,
\[
q+\sum_{N=1}^\infty \left[\prod_{i=0}^{N-2} \left(1+q^{2^i}+q^{-2^i}\right)\right] \left(1+q^{2^{N-1}}\right)q^{2^N}
\]
is the generating function of $\overline{f}_\mathrm{BSD}(n)$.  This is the desired generating series because the $N$-th term in the series represents those ways of writing a number in short BSD representation with leading term $1\cdot2^N$.  Observe that $(1+q^{2^i}+q^{-2^i})=(q^{0\cdot2^i}+q^{1\cdot2^i}+q^{T\cdot2^i})=((q^{2^i})^0+(q^{2^i})^1+(q^{2^i})^T)$, so we can think of the exponent on $q^{2^i}$ as recording the digit in the $2^i$ place of the short BSD representation.
Since we exclude $1T$ as a possibility for the leading two terms, we use $1+q^{2^{N-1}}$ instead of $1+q^{2^{N-1}}+q^{-2^{N-1}}$ for the allowed digits in the $2^{N-1}$ place.  The extra $q$ in front is needed because otherwise, if we instead began the sum with $N=0$, the $N=0$ term of the sum would be $q+q^{3/2}$, which includes $q^{3/2}$ unnecessarily.

We will prove two claims, from which the generating function identity \eqref{generating_function_identity} will follow.
\begin{enumerate}
    \item A finite analogue of \eqref{generating_function_identity}, namely
    \begin{align}\label{eqn}
q\prod_{n=0}^{M-1} \left( 1+ q^{2^n}+q^{2\cdot 2^n}\right) = q+\sum_{N=1}^M \left[\prod_{i=0}^{N-2} \left(1+q^{2^i}+q^{-2^i}\right)\right] \left(1+q^{2^{N-1}}\right)q^{2^N},
    \end{align}
    is true for $M\in\mathbb{N}$.
    \item There is a function $r(M)$ that goes to infinity with $M$ so that the first $r(M)$ coefficients in the above finite identity are unchanged if $M$ is replaced by infinity.  That is, each side of the identity converges formally.
\end{enumerate}

For the first claim, we proceed by induction. For $M=1$, we have on the left side
\[
q\prod_{n=0}^0 \left( 1+ q^{2^n}+q^{2\cdot 2^n}\right) = q+q^2+q^3,
\]
and on the right side, we have
\[
q+\sum_{N=1}^1 \left[\prod_{i=0}^{N-2} \left(1+q^{2^i}+q^{-2^i}\right)\right] \left(1+q^{2^{N-1}}\right)q^{2^N}= q+\left(1+q^{2^0}\right)q^{2^1}=q+q^2+q^3,
\]
so the base case holds.  We now assume that the identity holds for some fixed $M$ and show that it holds for $M+1$. On the right side of Equation \eqref{eqn}, we have
\begin{align*}
    &q+\sum_{N=1}^{M+1} \left[\prod_{i=0}^{N-2} \left(1+q^{2^i}+q^{-2^i}\right)\right] \left(1+q^{2^{N-1}}\right)q^{2^N}\\
    &\qquad =q+\sum_{N=1}^M \left[\prod_{i=0}^{N-2} \left(1+q^{2^i}+q^{-2^i}\right)\right] \left(1+q^{2^{N-1}}\right)q^{2^N}\\
    &\qquad\qquad + \left[\prod_{i=0}^{M-1} \left(1+q^{2^i}+q^{-2^i}\right)\right]\left(1+q^{2^{M}}\right)q^{2^{M+1}}\\
    &\qquad =q\prod_{n=0}^{M-1} \left( 1+ q^{2^n}+q^{2\cdot 2^n}\right)+ \left[\prod_{i=0}^{M-1} \left(1+q^{2^i}+q^{-2^i}\right)\right]\left(1+q^{2^{M}}\right)q^{2^{M+1}}.
    \end{align*}
In the second term, we now multiply by $1=q^{2^M-1}q^{-2^M+1}$, distribute $q^{2^M-1}=q^{1+2+\cdots+2^{M-1}}$ across the large product, and use the fact that $\big(1+q^{2^i}+q^{-2^i}\big)q^{2^i} = 1+q^{2^i}+q^{2\cdot 2^i}$ to obtain
    \begin{align*}
    &q\prod_{n=0}^{M-1} \left( 1+ q^{2^n}+q^{2\cdot 2^n}\right)+ \left[\prod_{i=0}^{M-1} \left(1+q^{2^i}+q^{-2^i}\right)\right]\left(1+q^{2^{M}}\right)q^{2^{M+1}}\\
        &\quad=q\prod_{n=0}^{M-1} \left( 1+ q^{2^n}+q^{2\cdot 2^n}\right) + \left[\prod_{i=0}^{M-1} \left(1+q^{2^i}+q^{2\cdot2^i}\right)\right]\left(1+q^{2^{M}}\right)q^{2^{M+1}-2^{M}+1}\\
    &\quad =q\left[ \prod_{n=0}^{M-1} \left( 1+ q^{2^n}+q^{2\cdot 2^n}\right)\right] \left( 1+ \left(1+q^{2^M}\right)q^{2^{M+1}-2^{M}}\right)\\
    &\quad=q\prod_{n=0}^{M} \left( 1+ q^{2^n}+q^{2\cdot 2^n}\right).
\end{align*}
Thus we have proven the first claim by induction.

For the second claim, we examine the tails of each side of Equation \eqref{generating_function_identity}. First, the product
\[
\prod_{n=M}^\infty \left( 1+ q^{2^n}+q^{2\cdot 2^n}\right)
\]
contains a summand of $1$, and then the next smallest exponent is $2^{M}$. Thus the finite product on the left side of Equation \eqref{eqn} and the infinite product on the left side of Equation \eqref{generating_function_identity} agree up to the $q^{2^{M}}$ term, and, in fact, they also agree on the $q^{2^{M}}$ term due to the factor of $q$ outside the product.  On the other side, in the tail
\[
\sum_{N=M+1}^\infty \left[\prod_{i=0}^{N-2} \left(1+q^{2^i}+q^{-2^i}\right)\right] \left(1+q^{2^{N-1}}\right)q^{2^N},
\]
the smallest power of $q$ is $q^{2^M+1}$, which appears in the $N=M+1$ summand from selecting the $q^{-2^i}$ term every time in the product and selecting $1$ instead of $q^{2^{N-1}}$. These choices give $q^{2^{M+1}-1-2-\cdots - 2^{M-1}} =q^{2^M+1}$. Thus the finite sum on the right side of Equation \eqref{eqn} and the infinite sum on the right side of Equation \eqref{generating_function_identity} agree up to the $q^{2^M}$ term. This completes the proof of the second claim and thus the theorem overall.\qed
\end{Proof}

\section{Connections to BSD representations with a specified number of digits}\label{section4}

We begin this section by recalling Monroe's results from \cite{ML}, as discussed previously in Section \ref{introduction}.  Theorem \ref{Monroe relation} below is the main result in \cite{ML} and states that the number of ways to write $n$ in BSD representation using $i$ digits is equal to the $\left(2^i-n\right)$-th term of the Stern sequence.

Given a non-negative integer $n$ and any $i\geq\lfloor \log_2(n)\rfloor+1$, recall that $f_{\text{BSD}}(n,i)$ denotes the number of ways to write $n$ in BSD representation with $i$ digits.  This is the number of ways to write $n$ in the form 
\[
n=\sum_{j=0}^{i-1} \epsilon_j 2^j, \qquad \epsilon_j\in\{0,1,T\}.
\]

\begin{theorem}[Theorem 2 of \cite{ML}]\label{Monroe relation}
If $0<n<2^i$, then $f_\mathrm{BSD}(n,i)=s(2^i-n)$. 
\end{theorem}

We restate Theorem \ref{main3} here for convenience and then prove it.

\begin{thm}
Let $n\ge 1$ be an integer and let $i=\lfloor \log_2 n\rfloor+1$ denote the number of digits in the standard binary representation of $n$. Then 
    \[
    f_\mathrm{BSD}(n,i+1)-f_\mathrm{BSD}(n,i)=\overline{f}_\mathrm{BSD}(n).
    \]
\end{thm}

\begin{proof}
Let $n$ be a positive integer and let $i=\lfloor \log_2 n \rfloor +1$.  Then $i$ is the number of digits in the standard binary representation of $n$, so $2^{i-1}\le n<2^i$.  We will first show that each of the representations counted by $\overline{f}_\mathrm{BSD}(n)$ must use either $i$ or $i+1$ digits.

Suppose a short BSD representation of $n$ uses $i+k+1$ digits, for some $k\ge 1$.  Then $n=\sum_{j=0}^{i+k}\epsilon_j2^j$, where each $\epsilon_j\in\{0,1,T\}$ and $\epsilon_{i+k}\neq0$.  Since $n$ is positive, we know $\epsilon_{i+k}\neq T$, and thus $\epsilon_{i+k}=1$.  Write $n=2^{i+k}+\mathcal{S}$, where $\mathcal{S}=\sum_{j=0}^{i+k-1} \epsilon_j 2^j$.  Then
\[
\mathcal{S}=n-2^{i+k}<2^i-2^{i+k}=2^i-2^{i+k-1}-2^{i+k-1}=-2^{i+k-1}-\left(2^{i+k-1}-2^i\right)\leq-2^{i+k-1}. 
\]
Since the representation of $n$ is assumed to be short, we have that $\epsilon_{i+k-1}\neq T$. So the smallest possible value of $\mathcal{S}$ is $\mathcal{S}=\left[0T\cdots T\right]_2$ with $i+k-1$ consecutive $T$'s. Therefore, $\mathcal{S}\ge -2^{i+k-1}+1$, but we must have $\mathcal{S}<-2^{i+k-1}$. This contradiction means that a short BSD representation of $n$ uses at most $i+1$ digits.

Suppose a short BSD representation of $n$ uses $i-k$ digits, for some $k\ge 1$. The largest BSD representation using $i-k$ digits contains a $1$ in every position and is $[1 1  \cdots 1]_2=2^{i-k}-1$. However, this is strictly less than $2^{i-1}$, and we know that $2^{i-1}\le n$. Thus it is impossible for a short BSD representation of $n$ to use $i-k$ digits for $k\ge 1$.  Combining this with the conclusion of the previous paragraph, we see that every short BSD representation of $n$ must use either $i$ or $i+1$ digits.

Thus all of the short BSD representations of $n$ are counted by $f_\mathrm{BSD}(n,i+1)$, but $f_\mathrm{BSD}(n,i+1)$ also counts some expansions of the form $n=\left[1T\cdots\right]_2$, which are forbidden.  How many of the expansions counted by $f_\mathrm{BSD}(n,i+1)$ are of this forbidden form?  If an expansion counted by $f_\mathrm{BSD}(n,i+1)$ is of the forbidden form, the expansion begins with $2^{i}-2^{i-1}$, which is equal to $2^{i-1}$.  Thus there is a corresponding expansion of $n$ starting with $0\cdot2^{i}+1\cdot 2^{i-1}$, and this expansion would be counted by $f_\mathrm{BSD}(n,i)$.  Hence the expansions counted by $f_\mathrm{BSD}(n,i+1)$ and of the forbidden form are in bijection with the expansions counted by $f_\mathrm{BSD}(n,i)$.  Thus $\overline{f}_\mathrm{BSD}(n)=f_\mathrm{BSD}(n,i+1)-f_\mathrm{BSD}(n,i)$.
\end{proof}

Observe that Theorem \ref{main3} leads to yet another proof of Theorem \ref{main2}, our statement that $\overline{f}_\mathrm{BSD}(n)=f_\mathrm{HB}(n-1)$ for all integers $n\geq 0$.

\begin{Proof}
By Corollary 2.1 in \cite{StolarskyDilcher}, 
\[
s(2^{i+j}- n)=s(2^i-n)+s(n)s(2^j- 1), \; \text{for }0\leq n\leq 2^i.
\]
Plugging in $j=1$, we get
\begin{align*}
s(2^{i+1}-n)&=s(2^i-n)+s(n)s(1)\\
&=s(2^i-n)+s(n).
\end{align*}
Rearranging, we obtain $s(n)=s(2^{i+1}-n)-s(2^i-n)$.  Using Theorem \ref{Monroe relation}, we have that $s(n)=f_\mathrm{BSD}(n, i+1)-f_\mathrm{BSD}(n,i)$, where the left side is equal to $f_\mathrm{HB}(n-1)$ and the right side is equal to $\overline{f}_\mathrm{BSD}(n)$ by Theorem \ref{main3}.
\end{Proof}

\end{document}